# An Enhanced Modelling Approach for Warehouse Sharing Platform System Designing Problem

Zeren XING [a], Yuehui WU [b,*], Shuangyuan YU [c]

a Department of Civil Engineering, Graduate School of Engineering, the University of Tokyo, Tokyo, 1138656, Japan

b College of Transportation Engineering, Dalian Maritime University, Dalian, 116026, China

c Department of Global Engineering, Undergraduate School of Engineering, Kyoto University, 6068317, Kyoto, Japan

**Abstract**

With the increasing importance of sustainability, warehouse sharing arises as a possible way to improve the efficiency of the existing logistics system. This paper studied the warehouse sharing platform systems (WSPS) and investigated its supply chain network, including factories, warehouses, and customers. We proposed an enhanced modelling approach for the warehouse sharing platform system design problem (WSPSDP) using the multi-allocation hub location routing problem framework. New elements such as inter-warehouse transportation and multiple-allocation scheme were added compared to the existing WSPS model. Then an adaptive large neighbourhood decomposition search heuristic was applied to solve our problem. Computational experiments were conducted on different-sized instances for comparison with the WSPS model without inter-warehouse transportation and the WSPS model with single-allocation scheme. The results suggested that our proposed WSPSDP model is more cost-efficient than the existing WSPS models, and it has the potential to promote the utilisation of existing cheap idle warehouses.

## 1. INTRODUCTION

With the increasing importance of sustainability, logistics has become an essential focus for governments from the future planning perspective (Raimbault et al., 2012). According to a report from the Japanese government, the changing business environment has led to dispersed logistics facilities, low average load factor, and small lot size (MLIT, 2020a), increasing urban freight traffic. The "Logistics Comprehensive Efficiency Law" (MLIT, 2020b) provided guidelines to resolve the inefficiency issues in logistics, wherein the efficiency of consolidated deliveries and supply chain cooperation was highlighted (Qureshi et al., 2022). Under this background, warehouse sharing platform systems (WSPS) arise as online platforms that connect manufacturers needing warehouse space and companies with excess warehouse capacity.

The conception of warehouse sharing is developed based on the principles of the sharing economy and platform businesses. Generally, warehouse sharing can reduce warehousing costs (labour and storage costs) and increase warehouse service levels, thereby improving supply chain efficiency. Moreover, it can also enhance the sustainability of the logistics systems (Ding & Kaminsky, 2020; Tornese et al., 2020). As a specific realisation of warehouse sharing conception, WSPS enables a dynamic supply network to respond to changing demand requirements. More significantly, WSPS can increase the usage of existing idle warehouses and relieve the inundation of manufacturer-owned warehouses (Pazour & Unnu, 2018; Qureshi et al., 2022).

This research focuses on designing the supply chain network of WSPS. The target network consists of factories, warehouses, and customers, where products are collected from the factories, temporarily stored in the warehouses, and finally delivered to the customers, as shown in Figure 1(a). Compared to the existing network, we include inter-warehouse transportation as it is usually cheaper than local vehicle routing due to the scale economics (de Camargo et al., 2013) so that the target network can be more cost-efficient. Moreover, we allow a multiple-allocation scheme between warehouses and factories/customers to improve the collection/delivery process further. Please find the comparison between the enhanced network and the existing one (the one in Qureshi et al. (2022)) in Figure 1. For notation simplicity, we name the enhanced network planning problem as the warehouse sharing platform system design problem (WSPSDP).

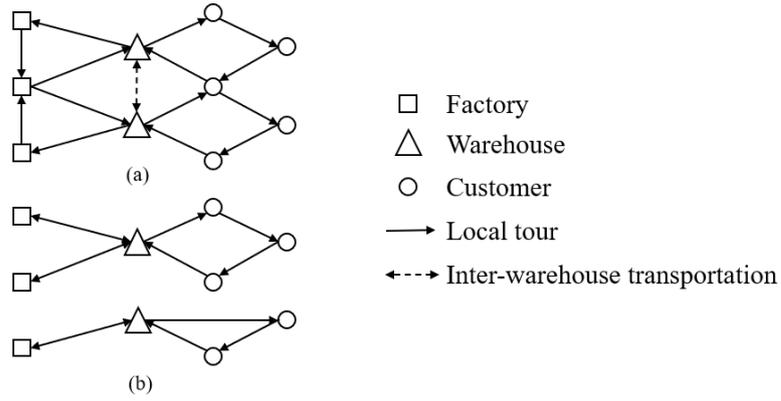

Figure 1. Enhanced network (a) and existing network (b)

The WSPSDP is a special case of a classical optimisation problem called the multi-allocation hub location routing problem (MAHLRP), which incorporates hub location and vehicle routing. For the details of the MAHLRP, please refer Fontes and Goncalves (2019) and Wu et al. (2022). The differences between the WSPSDP and the MAHLRP mainly lay in the following three points: i) The utilisation of warehouses incurs variable costs depending on the inbound flows rather than the fixed costs, as space in warehouses is rented for a short time; ii) The collection and delivery processes are distinct, i.e., each local tour should visit factories only or customers only; iii) One assumption is that each pair of factory/customer and warehouse is connected by at most one vehicle. Thus, in this study, we first formulate the WSPSDP as a special case of the MAHLRP and then introduce an adaptive large neighbourhood decomposition search (ALNDS) algorithm to solve it. Finally, we compare the WSPSDP and the existing WSPS network without inter-warehouse transportation, named WSPS-WI, as well as the WSPS network with single-allocation scheme (WSPS-SA) using different instances to demonstrate the efficiency and usefulness of the proposed model and algorithm.

The remainder of this paper is organised as follows: Section 2 reviews the existing literature about the WSPS and the hub location routing problem (HLRP); Section 3 proposes the mathematical model used in this research. The ALNDS algorithm used to solve the WSPSDP is introduced in Section 4. Section 5 describes the generated instances and settings of the computational experiments. The computation results are reported and discussed in Section 6. Finally, Section 7 gives the conclusions and future works of this research.

## 2. LITERATURE REVIEW

In the field of WSPS, little literature investigates its whole supply chain network. Pazour and Unnu (2018) presented a mathematical formulation for calculating the cost for a simple WSPS model without considering delivery processes to customers. In their study, only direct links were allowed within the network. The total cost was considered as the sum of fixed costs for setting up the warehouse, the variable costs for sorting products and transportation costs. They used a mixed integer programming approach to solve the problem and compared the costs to the case without on-demand warehouse sharing. More recently, Tornese et al. (2020) described the main features of WSPS and compared it with other warehousing types, representing a starting point of the WSPS research. They proposed a framework for the WSPS problem but did not provide a detailed mathematical model to formulate this framework. Later, Qureshi et al. (2022) proposed a mathematical model for the WSPS problem. They considered the WSPS as a two-level problem. In the upper level, an assignment problem between factories and warehouses was considered; in the lower level, a location routing problem between warehouses and customers was considered. They used a Tabu search heuristic to solve a series of scenarios based on hypothetical test instances. Moreover, they compared the results to the cases without splitting factory production and vehicle routings. The results showed that by considering the splitting of the factory production and vehicle routings for the outward delivery of the stored products in the warehouses, the total cost and $CO_2$ emissions can be significantly reduced.

In this study, we used the mathematical model of MAHLRP, which is a branch of HLRP, to formulate the WSPSDP. HLRP can be categorised into two types: Single-allocation hub location routing Problem (SAHLRP) and MAHLRP. Nagy and Salhi (1998) first proposed the HLRP with the name of many-to-many location-routing problem. They proposed an integer linear programming formulation for this problem and utilised a *locate first–route second* heuristic algorithm to solve it on a single instance with 249 non-hub nodes. de Camargo et al. (2013) introduced a new SAHLRP model with simultaneous collections and deliveries. They assumed that a fixed cost was imposed upon the hubs and vehicles. Moreover, they decomposed the problem into two sub-problems: a transportation problem and a feasibility problem. Then the problem was solved exactly by a tailored Benders decomposition algorithm. The results were compared to the CPLEX MILP solver, and they concluded that this method can solve instances up to 100 nodes. Karimi (2018) studied a capacitated SAHLRP with simultaneous pickup and delivery. The study introduced a polynomial-size mixed integer programming formulation for the problem including several valid inequalities. Moreover, a tabu-search based heuristic was proposed to solve the problem. The results from computational tests showed that the proposed valid inequalities and algorithm works well for their SAHLRP model. Yang et al. (2019) investigated the capacitated single-allocation hub location-routing problem (CSAHLRP) with distinct collection and delivery processes. Moreover, they proposed a new MILP model and developed a memetic algorithm (MA) to solve larger-sized problems. Their computational experiments showed that MA can find feasible solutions to the problem for all types of instances in a reasonable computing time and with limited gaps compared to the lower bounds of CPLEX. Yang (2018) considered the impacts of $CO_2$ emissions on the CSAHLRP. A bi-objective mathematical model was introduced to minimise the total costs and $CO_2$ emission at the same time. Then a bi-objective MA algorithm was used as the solution approach. The results showed that the method was efficient in obtaining high-quality solutions for small instances in reasonable time.

The amount of research conducted on the MAHLRP is relatively limited compared with the SAHLRP. Çetiner et al. (2010) conducted the initial investigation on the MAHLRP. Their contribution was restricted to the proposal of the MAHLRP framework without presenting a mathematical formulation of the problem. In their study, both hubs and vehicles were assumed to be uncapacitated. Moreover, they developed an iterative heuristic algorithm. The first step involves identifying hub locations, and multiply allocating postal offices to the hubs. In the second step, the offices were routed via local tours originating from the allocated hubs. Finally, a case study was performed using the Turkish postal delivery system data. Basirati et al. (2020) investigated a MAHLRP with simultaneous collection and delivery process within hard time window. They considered that hubs and fleets have limited capacity. Furthermore, they proposed a bi-objective model to balance travel costs among different routes and to minimise the total sum of fixed costs of locating hubs, the costs of handling, travelling, assigning, and transportation costs. Small-and-medium-sized problems were solved using an augmented $\varepsilon$-constraint technique; on the other hand, large size problems were solved utilising the proposed multi-objective imperialist competitive algorithm (MOICA). The results suggested that MOICA performs better than the non-dominated sorting genetic algorithm (NSGA-II) for large-sized problems. Wu et al. (2022) conducted an application of MAHLRP for the design of an intra-city express service system. They proposed a mixed integer programming mathematical model with simultaneous collection and delivery processes and utilised an ALNDS algorithm to solve the problem. Furthermore, computational experiments were conducted using Australian Post data to compare the proposed MAHLRP model with the SAHLRP; simultaneously, the ALNDS heuristic algorithm was compared to the CPLEX solver. The results showed that the MAHLRP can efficiently reduce the operating cost compared to the SAHLRP, and their ALNDS algorithm performed better than CPLEX on both solution quality and computational time.

In conclusion, there is a lack of literature related to WSPS; meanwhile, the existing WSPS model is considered inefficient and can be improved from several perspectives. Our study proposes and includes some new aspects to develop a new WSPSDP model. In addition, it is the first to formulate the WSPS model based on the MAHLRP, which is one of our contributions to the literature.

## 3. PROBLEM DESCRIPTION

In this research, the WSPSDP problem is defined on a complete directed graph $G = (N, A)$ with node set $N = F \cup C \cup W$ and arc set $A$. Here, $F$, $C$ and $W$ represent factory set, customer set and warehouse set, respectively. The flow between factories $i \in F$ and customers $j \in C$ is denoted as $q_{ij}$, and it is transported via warehouses.

For each arc $(i,j) \in A | i,j \in N$, $c_{ij}$ represents its travel distance. The local tour cost depends on the travel distance only. On the other hand, both travel distance and transferred flow determine the inter-warehouse transportation cost (Wu et al., 2023; Yang et al., 2019). $\alpha$ and $\beta$ are the unit costs for the inter-warehouse transportation and local tour, respectively. Inter-warehouse transfer is assumed not subject to capacity restrictions, whose economies of scale can be reflected by $\alpha$. Vehicle capacity is denoted by $Q^v$, and warehouse capacity is denoted by $Q^m$ ($m \in W$).

The unit warehouse variable cost is denoted by $a_m$, which consists of warehouse rental costs and freight operation costs. Moreover, it is assumed that the variable cost only occurs for the original warehouse, i.e., if products are transferred through warehouses $m$ and $n$, then only the variable cost for warehouse $m$ is considered.

In addition, for each arc $(i,j) \in A | i,j \in N$ and warehouse $m \in W$, $x_{ijm}$ denotes a binary variable which takes the value of 1 if there is a vehicle leaving from warehouse $m$ and travelling directly from node $i$ to node $j$. $z_{ijmn}$ indicates the fraction of flow from factory $i \in F$ to customer $j \in C$ traversing through warehouses $m \in W$ and $n \in W$. $u_{im}$ indicates the collection load on the vehicle leaving from warehouse $m \in W$ after visiting node $i$; and $v_{im}$ indicates the delivery load on it just before visiting node $i$. $M$ is a large positive number. Using the above parameters and assumptions, the WSPSDP is defined as follows:

$$\min \sum_{i \in F} \sum_{j \in C} \sum_{m \in W} \sum_{n \in W} a_m z_{ijmn} q_{ij} + \sum_{i \in N} \sum_{j \in N} \sum_{m \in W} \beta c_{ij} x_{ijm} + \sum_{i \in F} \sum_{j \in C} \sum_{m \in W} \sum_{n \in W} \alpha c_{mn} z_{ijmn} q_{ij} \quad (1)$$

$$s.t. \sum_{m \in W} \sum_{n \in W} z_{ijmn} = 1 \quad \forall i \in F, j \in C \quad (2)$$

$$\sum_{i \in F} \sum_{j \in C} \sum_{n \in W} z_{ijmn} q_{ij} + \sum_{i \in F} \sum_{j \in C} \sum_{n \neq m \in W} z_{ijnm} q_{ij} \leq Q^m \quad \forall m \in W \quad (3)$$

$$\sum_{j \neq i \in N} \sum_{m \in W} x_{ijm} \geq 1 \quad \forall i \in F \cup C \quad (4)$$

$$\sum_{j \neq i \in N} x_{ijm} \leq 1 \quad \forall i \in F \cup C, m \in W \quad (5)$$

$$\sum_{j \in N} x_{ijm} = \sum_{j \in N} x_{jim} \quad \forall i \in N, m \in W \quad (6)$$

$$M \sum_{j \neq i \in N} x_{ijm} \geq \sum_{j \in C} \sum_{n \in W} z_{ijmn} \quad \forall i \in F, m \in W \quad (7)$$

$$M \sum_{j \neq i \in N} x_{ijm} \geq \sum_{j \in F} \sum_{n \in W} z_{jinm} \quad \forall i \in C, m \in W \quad (8)$$

$$\sum_{i \in F} \sum_{j \in C} x_{ijm} + \sum_{i \in F} \sum_{j \in C} x_{jim} = 0 \quad \forall m \in W \quad (9)$$

$$\sum_{n \in W} \sum_{h \in W} x_{nhm} = 0 \quad \forall m \in W \quad (10)$$

$$u_{im} + \sum_{s \in C} \sum_{n \in W} q_{js} z_{jsmn} - Q^v(1 - x_{ijm}) \leq u_{jm} \quad \forall i \in N, j \neq i \in F, m \in W \quad (11)$$

$$v_{im} - \sum_{s \in F} \sum_{n \in W} q_{si} z_{sinm} + Q^v(1 - x_{ijm}) \leq v_{jm} \quad \forall i \in C, j \neq i \in N, m \in W \quad (12)$$

$$v_{im} + u_{im} \leq Q^v \ \forall i \in F \cup C, m \in W \tag{13}$$

$$x_{ijm} = \{0,1\} \ \forall i,j \in N, m \in W \tag{14}$$

$$z_{ijmn} \geq 0 \ \forall i \in F, j \in C, m \in W, n \in W \tag{15}$$

$$u_{im} \geq 0; \ v_{im} \geq 0 \ \forall i \in N, m \in W \tag{16}$$

In the above model, objective function (1) minimises the supply chain network's total cost, consisting of three terms: warehouse variable cost, local tour cost, and inter-warehouse transportation cost. Constraint (2) ensures that all flows between factories and customers are served. Constraint (3) limits warehouse capacity. Constraint (4) guarantees that each factory and customer should be visited at least once. Constraint (5) ensures that each pair of factory/customer and warehouse should be linked by at most one vehicle. Constraint (6) is the vehicle flow conservation constraint. Constraints (7) and (8) indicate that flow should occur only between linked factories/customers and warehouses. Each route should visit factories only or customers only, which is guaranteed by Constraint (9). Constraint (10) imposes that local tours cannot link warehouses. Constraints (11) and (12) specify the collection and delivery load on vehicles, respectively. Constraint (13) impose limitations on the vehicle's capacity. Constraints (14)-(16) are variable domains.

## 4. ADAPTIVE LARGE NEIGHBOURHOOD DECOMPOSITION SEARCH

T The WSPSDP consists of three subproblems: warehouse selection, non-warehouse node allocation, and vehicle routing. As a result, its solution space is very huge. Therefore, the ALNDS meta-heuristics algorithm proposed by Wu et al. (2022) is utilised to solve our problem. The flowchart of the algorithm is illustrated in Figure 2.

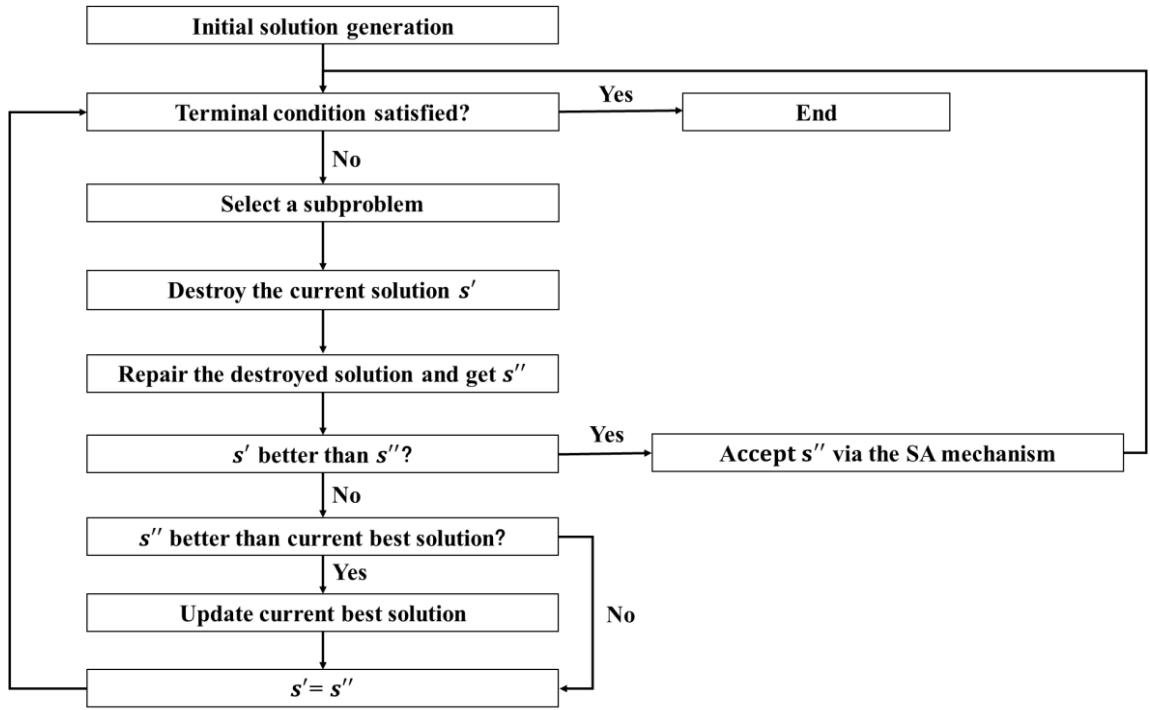

Figure 2. Flowchart of the introduced algorithm

The ALNDS first builds an initial solution and iteratively improves it using the destroy and repair operators. The initial solution is constructed in two phases; in the first phase, a greedy algorithm is applied to assign the customers and factories to warehouses, whereas a nearest neighbour algorithm is used to construct routes in the second phase, as shown in Algorithm 1.

In each iteration, the algorithm selects one of the subproblems, and the neighbourhood of the current solution is explored using destroy and repair operators. For each subproblem, a series of destroy and repair operations is employed. For example, the warehouse location subproblem uses random removal and random opening operators. To incorporate the multi-allocation scheme, copies of non-warehouse nodes (sub-nodes) are created and added/removed using duplication/ deduplication operators. The vehicle routes are also destroyed by removing subnodes based on random selection, by selecting worst-cost customers giving the highest removal gains, or by removing similar sub-nodes (Shaw removal (Ropke & Pisinger (2006))) and so on. Please find more details of the used operators in Wu et al. (2022).

> **Algorithm 1:** Two-phase construction heuristic algorithm
> **While** unassigned non-warehouse nodes remained:
>     Select unassigned non-warehouse node $i$ with maximum demand
>     **If** $i$ is a factory
>         Assign it to the cheapest warehouse that has enough capacity
>         Update the remaining capacity of the related warehouse
>     **Else**
>         Assign it to the nearest warehouse
>     **Endif**
> **End While**
> **For** each warehouse $m$:
>     **While** unlinked factories assigned to warehouse $m$ remained:
>         Set current node $v = m$
>         Set current vehicle capacity $q = Q^v$
>         **While** feasible unlinked factories remained:
>             Select unlinked factories $i$ that is nearest to $v$
>             $v = i$
>             $q = q - \sum_{j \in C} q_{ij}$
>         **End While**
>     **End While**
>     **While** unlinked customers assigned to warehouse $m$ remained:
>         Set current node $v = m$
>         Set current vehicle capacity $q = Q^v$
>         **While** feasible unlinked customers remained:
>             Select unlinked customers $i$ that is nearest to $v$
>             $v = i$
>             $q = q - \sum_{j \in C} q_{ji}$
>         **End While**
>     **End While**
> **End For**
> **Output** generated initial solution

Each subproblem, destroy operator and repair operator is addressed with a weight and selected by a roulette mechanism according to their weights. For example, the selection possibility of destroy operator $o$ is calculated by Eq. (17), in which $w_o$ is its weight, and $O$ is the destroy operator set. Note that the subproblem and repair operators are selected independently of the destroy operators and vice versa.

$$P(o) = \frac{w_o}{\sum_{k \in O} w_k} \tag{17}$$

The weights are adjusted based on their previous performance, and we follow the adaptive weight control strategy introduced by Ropke and Pisinger (2006). The entire search is divided into a number of segments, and a *segment* is a number of iterations of the ALNS heuristic in which the weights remain fixed. At the end of each segment, the weights are adjusted based on the performance. To evaluate the performance of the operators, we address each operator with a score, and these scores are set as zero at the beginning of the segment. Whenever an operator is applied, its score is increased by either $\sigma_1, \sigma_2, \sigma_3$ or $0$ ($\sigma_1 > \sigma_2 > \sigma_3 > 0$) according to the quality of the solution obtained via it. For instance, let $w_{o,seg}$ and $w_{o,seg+1}$ denote the weight of destroy operator $o$ at the beginning and end of the segment, respectively. $w_{o,seg+1}$ is obtained from $w_{o,seg}$ as follows:

$$w_{o,seg+1} = (1-\eta)w_{o,seg} + \eta \frac{score_o}{time_o} \qquad (18)$$

where $score_o$ and $time_o$ represent the score and selected times of operator $o$ in the segment. $\eta$ is a parameter range from [0,1] and is used to reflect the impact of previous weights on the newly calculated weights.

When a new-obtained solution is better than the existing one, it is adopted as the input for the subsequent iteration. Otherwise, to decide whether to accept it or not, a simulated annealing mechanism is implemented. The probability of a degrading solution $s'$ ($f(s') > f(s_{current})$) is calculated as below:

$$p = e^{-(f(s')-f(s_{current}))/T} \qquad (19)$$

where $T$ is called as temperature and set as $T_{start}$ when the algorithm starts to ensure that a solution 20% worse than the initial solution is accepted with a possibility of 30%. In each iteration, $T$ is reduced by multiplying a cooling factor $\theta \in \{0,1\}$ to maintain the stability of the search. The search terminates when the given number of iterations is reached. For more details, readers are referred to Wu et al. (2022).

## 5. COMPUTATIONAL EXPERIMENTS

### 5.1 Instance Generation

We evaluated the efficiency of the proposed WSPSDP model on the instances generated from a public data set, namely Turkish network (https://ie.bilkent.edu.tr/~bkara/dataset.php). The Turkish network contains data on 81 cities in Turkey. We directly use the distances and flow between the nodes (or demand centres). Furthermore, to be consistent with the literature, we applied the same 16 candidate warehouse locations as Yaman et al. (2007) and Alumur et al. (2009).

We created several test instances based on the Turkish network data introduced above. Small-sized instances have 18 instances with up to 20 non-warehouse nodes and 5 or 7 warehouses. Regarding medium-and-large-sized instances, they contain 20 instances ranging from 25 to 50 non-warehouse nodes and 7 warehouses. We selected the warehouse candidates as below:

(i) For instances with 7 warehouses, the warehouse candidates were randomly selected from the above referred 16 candidate warehouse locations, while factory and customer nodes were randomly selected from the remaining 65 nodes.

(ii) For instances with 5 warehouses, the warehouse candidates were randomly selected from the 7 selected warehouses of the corresponding instances.

Figure 3 shows an instance with 7 warehouses, 5 factories, and 25 customers. Here, red rectangles, green triangles, and blue circles represent warehouses, factories, and customers.

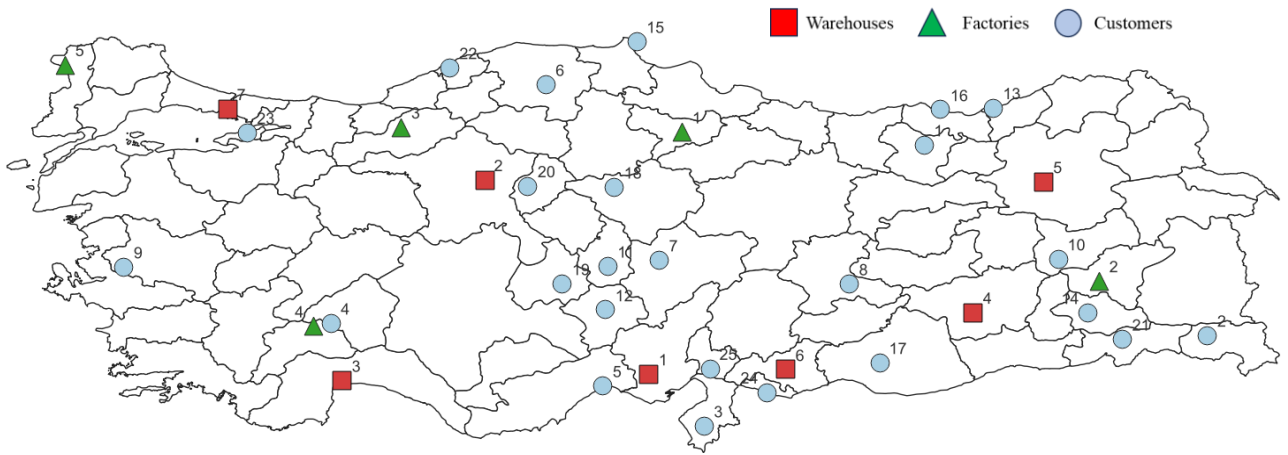

Figure 3. Example instance with 7 warehouses, 5 factories, and 25 customers

Moreover, the Turkish network data does not include the data on variable costs, warehouse capacities and vehicle capacities. Therefore, we randomly generated the variable costs for all candidate warehouses from 0.1 to 0.3; and performed computational experiments under different warehouse capacity conditions. For small-sized instances, three types of warehouse capacities were considered, including 'Compact (C)', 'Moderate (M)' and 'Spacious (S)', which corresponds to 0.3, 0.5 and 0.7 times of the total demand, respectively. Medium-and-large-sized instances contain four types of warehouse capacities: 'Compact (C)', 'Moderate (M)', 'Spacious (S)', and 'Large (L)', which corresponds to 0.3, 0.5, 0.7 and 2 times of the total demand, respectively.

We named our instances as *W-F-C-$Q^m$*, where *W* indicates the number of candidate warehouses, *F* stands for the number of factories, *C* represents the number of customers, and $Q^m$ indicates the type of warehouse capacities. On the other hand, regarding the vehicle capacity values, we referred to the study by Karimi (2018), where the vehicle capacities were determined by taking a random value between the maximum demand among all the non-warehouse nodes and the total demand. Furthermore, unit inter-warehouse transportation cost $\alpha$ was assumed to be 0.001, and the unit local tour cost $\beta$ was set as 1 to balance the three segments of the objective function, ensuring that no single term is significantly larger than the others.

To examine whether the proposed model performs well by adding the multiple-allocation scheme and inter-warehouse transportation, we compared our proposed WSPSDP model to the WSPS-SA and the WSPS-WI for each instance.

**5.2 Results and Discussions**

The program code for the solution algorithm was written using JAVA, and the computational experiments were performed on a PC with an Intel i7-10750H CPU and a memory of 64 GB. Furthermore, CPLEX 22.1.1.1 was utilised to solve the proposed WSPSDP model for the small instance in up to 10800 seconds. For each instance, the ALNDS heuristic was executed 10 times. The results were obtained and compared to the WSPS-SA and WSPS-WI models.

**5.2.1 Small-sized instances**

In this sub-section, we present and discuss the results for small-sized instances. Table 1 shows the results for small-sized instances. The notations used in Table 1 are defined as follows:

$S_{best}$: The best objective value of 10 executions by the ALNDS heuristic for the WSPSDP model.

$S_{ave}$: The average objective value of 10 executions by the ALNDS heuristic for the WSPSDP model.

$T_M$ (s): CPU time of ALNDS heuristic in seconds to obtain the best objective value.

$N_M$: The number of multi-allocation non-warehouse nodes.

$UB$: The best objective value obtained by CPLEX in 10800 seconds.

$T_C$ (s): CPU time of CPLEX in seconds to obtain the best objective value.

$\%S_{UB}$: Difference between $S_{best}$ and $UB$, divided by $UB$ (in %), that is, $\%S_{UB} = (UB - S_{best})/UB \times 100\%$.

$S_{best}^{SA}$: The best objective value of 10 executions by the ALNDS heuristic for the WSPS-SA case.

$\%S_{SA}$: Difference between $S_{best}$ and $S_{best}^{SA}$, divided by $S_{best}^{SA}$ (in %), that is, $\%S_{SA} = (S_{best}^{SA} - S_{best})/S_{best}^{SA} \times 100\%$.

$S_{best}^{WI}$: The best objective value of 10 executions by the ALNDS heuristic for the WSPS-WI case.

$\%S_{WI}$: Difference between $S_{best}$ and $S_{best}^{WI}$, divided by $S_{best}^{WI}$ (in %), that is, $\%S_{WI} = (S_{best}^{WI} - S_{best})/S_{best}^{WI} \times 100\%$.

Table 1. Results for small-sized instance tests

| Instance | WSPSDP | | | | CPLEX | | | WSPS-SA | | WSPS-WI | |
|---|---|---|---|---|---|---|---|---|---|---|---|
| $W$-$F$-$C$-$Q^m$ | $S_{best}$ | $S_{ave}$ | $T_M$ (s) | $N_M$ | $UB$ | $T_C$ (s) | $\%S_{UB}$ | $S_{best}^{SA}$ | $\%S_{SA}$ | $S_{best}^{WI}$ | $\%S_{WI}$ |
| 5-5-5-C | 23597.16 | 23597.16 | 52.96 | 2 | 23597.16 | 106.34 | 0 | 23867.31 | 1.13 | 30031.14 | 21.42 |
| 5-5-5-M | 21817.78 | 21817.78 | 41.57 | 1 | 21817.78 | 108.78 | 0 | 21896.35 | 0.36 | 26427.00 | 17.44 |
| 5-5-5-S | 21447.98 | 21447.98 | 1.42 | 0 | 21447.98 | 61.58 | 0 | 21447.98 | 0 | 25020.90 | 14.28 |
| 7-5-5-C | 21708.59 | 21708.58 | 1.28 | 0 | 21708.59 | 36.25 | 0 | 21708.59 | 0 | 30036.97 | 27.73 |
| 7-5-5-M | 21520.63 | 21520.63 | 1.34 | 0 | 21520.63 | 72.30 | 0 | 21520.63 | 0 | 26251.00 | 18.02 |
| 7-5-5-S | 20844.18 | 20844.18 | 1.62 | 0 | 20844.18 | 17.12 | 0 | 20844.18 | 0 | 25225.31 | 17.37 |
| 5-5-10-C | 43450.65 | 43587.97 | 288.34 | 2 | 43450.65 | 2101.34 | 0 | 44856.04 | 3.13 | 56124.88 | 22.58 |
| 5-5-10-M | 40649.17 | 40649.17 | 328.49 | 1 | 40649.17 | 1446.70 | 0 | 41297.46 | 1.57 | 47352.03 | 14.16 |
| 5-5-10-S | 40474.81 | 40539.11 | 207.33 | 1 | 40474.81 | 2558.91 | 0 | 40618.75 | 0.35 | 44861.34 | 9.78 |
| 7-5-10-C | 24398.95 | 24513.01 | 486.25 | 1 | 25848.92 | 605.09 | 5.61 | 25969.23 | 6.05 | 36865.96 | 33.82 |
| 7-5-10-M | 22959.68 | 22959.68 | 322.85 | 1 | 22959.68 | 946.77 | 0 | 22998.57 | 0.17 | 32206.77 | 28.71 |
| 7-5-10-S | 22340.18 | 22361.77 | 328.28 | 1 | 25544.24 | 1163.89 | 12.54 | 22426.54 | 0.39 | 28493.89 | 21.60 |
| 5-5-15-C | 67516.01 | 67642.69 | 477.16 | 2 | 67527.41 | 10800 | 0.02 | 70129.73 | 3.73 | 82179.79 | 17.84 |
| 5-5-15-M | 62502.05 | 62645.95 | 452.59 | 1 | 62522.36 | 10800 | 0.03 | 63638.51 | 1.79 | 69867.19 | 10.54 |
| 5-5-15-S | 61568.12 | 61592.07 | 299.24 | 1 | 61568.12 | 10800 | 0 | 61696.36 | 0.21 | 65970.03 | 6.67 |
| 7-5-15-C | 60704.47 | 60811.52 | 754.45 | 2 | 61188.55 | 10800 | 0.79 | 62792.67 | 3.33 | 76341.99 | 20.48 |
| 7-5-15-M | 56546.62 | 56783.91 | 733.51 | 1 | 56986.60 | 10800 | 0.77 | 57391.77 | 1.47 | 65121.96 | 13.17 |
| 7-5-15-S | 54836.65 | 55096.43 | 761.73 | 1 | 54847.65 | 10800 | 0.02 | 55328.42 | 0.89 | 62557.07 | 12.34 |

From Table 1, we can observe that for 11 instances out of 18 in total, CPLEX can obtain optimal solutions within 10800 seconds. However, for all small-sized instance tests, the value of $\%S_{UB}$ ranges from 0 up to 12.54%. Moreover, Figure 4 (a) and (b) compare the average CPU time of the ALNDS heuristic and CPLEX with candidate warehouse numbers 5 and 7, respectively. The CPU time of CPLEX increased drastically as the instance size became larger, and it reached 10800 seconds when the instance size was equal to or greater than 20 non-warehouse nodes, regardless of the candidate warehouse number. The above results suggest that the ALNDS heuristic can find optimal or better solutions for the WSPSDP model within shorter CPU times.

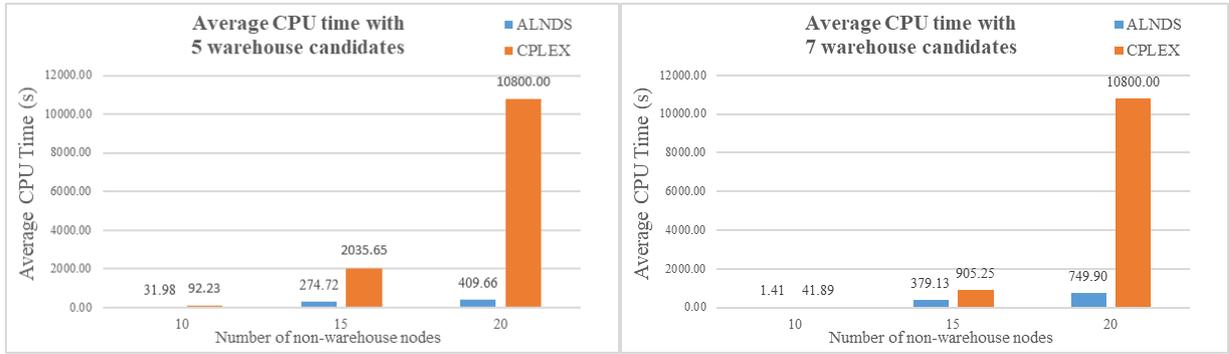

Figure 4. Average CPU time of ALNDS heuristic and CPLEX of

(a) Instances with 5 warehouse candidates (b) Instances with 7 warehouse candidates

Moreover, Table 1 demonstrates that for all small-sized instances, the values of $\%S_{WI}$ are positive, ranging from 6.67% up to 33.82%. The results indicate that for small-sized instances, our proposed WSPSDP model improved the supply chain efficiency of the WSPS network by adding inter-warehouse links, regardless of the candidate warehouse number. On the other hand, the values of $\%S_{SA}$ are non-negative, ranging up to 6.05%, suggesting that for small-sized instances, our proposed WSPSDP model overtakes the WSPS-SA network in terms of efficiency by adding the multi-allocation scheme. However, there are four instances with 5 factories and 5 customers, for which the values of $\%S_{SA}$ and $N_M$ are zeros. This means that WSPSDP obtains the same solution under particular conditions as WSPS-SA (no multi-allocation nodes exist).

The possible reasons are discussed by comparing the variable costs and geographic locations of the selected nodes. Table 2 presents the unit variable costs for instances '5-5-5-M' ($N_M = 1$) and '7-5-5-M' ($N_M = 0$). In addition, the vehicle routings of the best solutions for instances '7-5-5-M' and '5-5-5-M' were plotted in Figure 5 (a) and (b), respectively. In both figures, black solid lines represent local tours, and blue dash lines represent inter-warehouse transportation routes. In Figure 5 (b), red solid lines were used for the vehicle route passing the multi-allocation node (Factory 3). We can observe that regarding 'Factory 5', in instance '5-5-5-M', it uses 'Warehouse 3' with a unit variable cost of 0.149; however, in instance '7-5-5-M', it uses 'Warehouse 7' with a unit variable cost of 0.184. This indicates that the addition of inter-warehouse transportation in the WSPS network can make it feasible for non-warehouse nodes to use a closer warehouse instead of using a far-distanced warehouse, even if its variable cost is cheaper.

Table 2. Unit variable costs for the candidate warehouses in the selected instances

|  | Warehouse 1 | Warehouse 2 | Warehouse 3 | Warehouse 4 | Warehouse 5 | Warehouse 6 | Warehouse 7 |
| --- | --- | --- | --- | --- | --- | --- | --- |
| 5-5-5-M | 0.217 | 0.191 | 0.149 | 0.243 | 0.239 | - | - |
| 7-5-5-M | 0.217 | 0.191 | 0.149 | 0.243 | 0.239 | 0.288 | 0.184 |

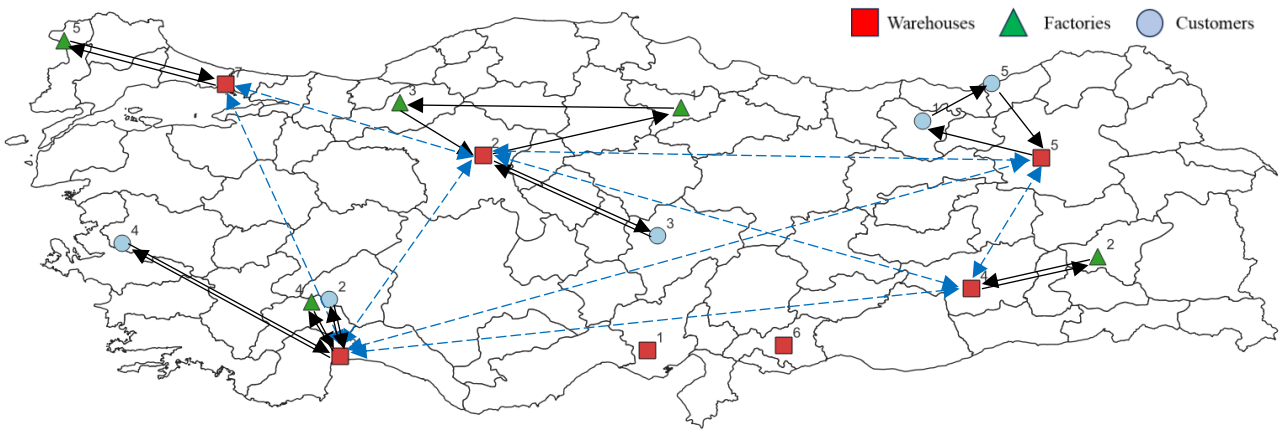

(a) Instance 7-5-5-M

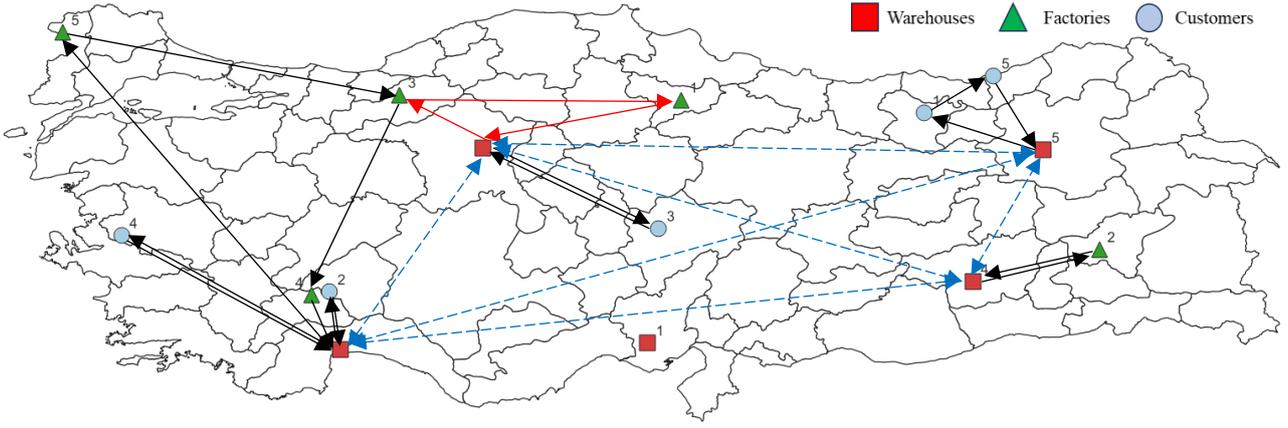

(b) Instance 5-5-5-M

Figure 5. Vehicle routing plots for the best solutions

### 5.2.2 Medium-and-large-sized instances

In this sub-section, we present and discuss the results for medium-and-large-sized instances. Table 3 presents the results for medium-and-large-sized instances. The notations used in Table 3 are the same as those in Table 1 except for the follows:

Used Warehouses: Numbers of warehouses used in $S_{best}$.

$\%R_{SD}$: Relative standard deviation of 10 times ALNDS heuristic executions for the WSPSDP case (in %), that is, $R_{SD} = SD/S_{ave} \times 100\%$. Here, $SD$ represents the standard deviation for the obtained 10 objective values. This indicator measures the deviation of the objective values of 10 executions disseminated around the mean value.

In Table 3, the values of $\%R_{SD}$ for all the test instances are less than 1%, and the average of $\%R_{SD}$ is 0.22% across all the instances. This proves that the ALNDS heuristic is stable and robust for providing high-quality solutions for medium-and-large-sized instances.

Table 3 demonstrates that for all medium-and-large-sized instances, the values of $\%S_{WI}$ are positive, ranging from 0.07% up to 25.30%. The results prove that our proposed WSPSDP model outperforms the existing WSPS-WI network by adding inter-warehouse links for medium-and-large-sized instances. On the other hand, the values of $\%S_{SA}$ are non-negative, ranging up to 7.93%. The results suggest that for medium-and-large-sized instances, our proposed WSPSDP model overtakes the WSPS-SA network in terms of efficiency by adding the multi-allocation scheme. However, there are five instances in which the WSPSDP obtains the same solution as the WSPS-SA when the warehouse capacities were set to 'Large'. Further analyses of the warehouse capacity effects are discussed in the following sub-section.

Table 3. Results for medium-and-large-sized instances tests

| Instance $W$-$F$-$C$-$Q^m$ | WSPSDP | | | | | | WSPSDP-SA | | WSPSDP-WI | |
|---|---|---|---|---|---|---|---|---|---|---|
| | $S_{best}$ | Used Warehouses | $T_M$ (s) | $S_{ave}$ | %$R_{SD}$ | $N_{MA}$ | $S_{best}^{SA}$ | %$S_{SA}$ | $S_{best}^{WI}$ | %$S_{WI}$ |
| 7-5-20-C | 52219.98 | 1,2,3,5,6,7 | 1053.12 | 52257.46 | 0.16 | 1 | 52407.38 | 0.36 | 69905.45 | 25.30 |
| 7-5-20-M | 51244.16 | 1,2,3,5,6,7 | 935.18 | 51336.16 | 0.22 | 1 | 51474.17 | 0.45 | 66232.62 | 22.63 |
| 7-5-20-S | 50730.55 | 1,2,3,4,5,6,7 | 932.36 | 50789.70 | 0.15 | 1 | 51162.67 | 0.84 | 51376.48 | 1.26 |
| 7-5-20-L | 51162.67 | 1,2,3,5,6,7 | 909.64 | 51162.67 | 0 | 0 | 51162.67 | 0 | 51376.48 | 0.42 |
| 7-5-25-C | 92919.71 | 1,2,3,4,5,6,7 | 1651.17 | 93654.68 | 0.49 | 3 | 97355.75 | 4.56 | 124358.28 | 25.28 |
| 7-5-25-M | 86016.03 | 1,2,3,5,6,7 | 1780.63 | 86081.73 | 0.09 | 1 | 87536.83 | 1.74 | 103284.54 | 16.72 |
| 7-5-25-S | 84463.76 | 1,2,3,5,6,7 | 1430.51 | 84523.67 | 0.11 | 1 | 84890.37 | 0.50 | 98537.86 | 14.28 |
| 7-5-25-L | 81951.88 | 1,2,3,5,6,7 | 1511.98 | 81951.88 | 0 | 0 | 81951.88 | 0 | 83782.23 | 2.18 |
| 7-5-30-C | 80067.89 | 1,3,4,5,6 | 2283.62 | 80433.84 | 0.32 | 1 | 86963.08 | 7.93 | 103133.01 | 22.36 |
| 7-5-30-M | 65685.90 | 1,5,6 | 2259.47 | 65914.73 | 0.43 | 1 | 70576.25 | 6.93 | 81982.78 | 19.88 |
| 7-5-30-S | 64752.95 | 1,5,6 | 1541.88 | 64990.94 | 0.25 | 1 | 65227.84 | 0.73 | 70411.56 | 8.04 |
| 7-5-30-L | 63627.70 | 1,6 | 1559.96 | 63668.74 | 0.05 | 0 | 63627.70 | 0 | 64855.74 | 1.89 |
| 7-10-30-C | 122097.04 | 1,2,5,6,7 | 2564.78 | 122559.31 | 0.24 | 2 | 122993.52 | 0.73 | 152385.45 | 19.88 |
| 7-10-30-M | 118622.42 | 1,2,5,6,7 | 1759.88 | 118857.17 | 0.21 | 1 | 119352.02 | 0.61 | 145822.46 | 18.65 |
| 7-10-30-S | 117079.30 | 1,2,5,6,7 | 2622.26 | 117369.14 | 0.22 | 1 | 118078.89 | 0.85 | 141192.14 | 17.08 |
| 7-10-30-L | 115839.12 | 6,7 | 1979.08 | 116213.85 | 0.29 | 0 | 115839.12 | 0 | 116026.73 | 0.16 |
| 7-10-40-C | 154438.76 | 1,2,5,6,7 | 3802.65 | 154863.21 | 0.30 | 2 | 157112.45 | 1.70 | 190648.09 | 18.99 |
| 7-10-40-M | 149566.19 | 1,2,5,6,7 | 4774.39 | 149985.08 | 0.29 | 2 | 150845.67 | 0.85 | 170736.19 | 12.40 |
| 7-10-40-S | 148806.93 | 1,2,3,5,6,7 | 4403.21 | 149076.78 | 0.26 | 1 | 150106.13 | 0.87 | 163804.55 | 9.16 |
| 7-10-40-L | 145642.76 | 6,7 | 3352.64 | 146036.58 | 0.35 | 0 | 145642.76 | 0 | 145740.93 | 0.07 |

**5.2.3 Effects of warehouse capacities**

In this sub-section, we discuss the effects of warehouse capacities on the solution of the WSPSDP. As mentioned in the previous sub-section, there are five instances that the WSPSDP yields the same solution as the WSPS-SA when the warehouse capacities were set to 'Large'. Example plots for the vehicle routings of the best solutions for instances '7-5-25-C' ($N_M = 3$) and '7-5-25-L' ($N_M = 0$) are shown in Figure 6 (a) and (b), respectively. The meanings of lines are the same as in previous figures. Table 4 presents the unit variable costs for the candidate warehouses of both instances.

Comparing the vehicle routings for Figure 6 (a) and (b), we observe that, for instance '7-5-25-L', only 'Warehouse 7' with the lowest unit variable cost was selected for the collection processes, and the closest warehouses were used for delivery processes to customers via inter-warehouse transportation. The warehouse with the lowest unit variable cost can be selected without considering the capacity constraint. Therefore, under 'Large' warehouse capacity conditions, the WSPSDP model solutions were the same as those from the WSPS-SA model.

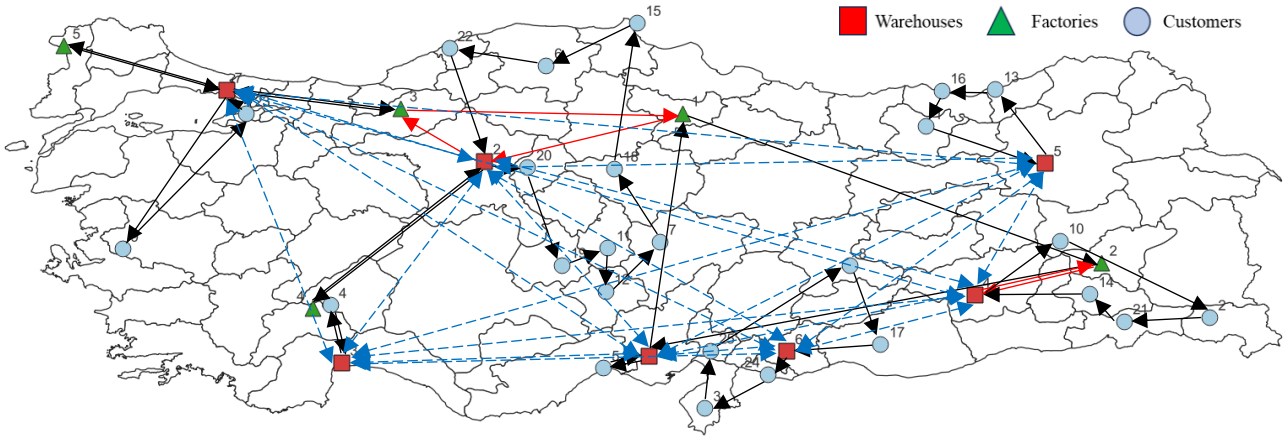

(a) Instance 7-5-25-C

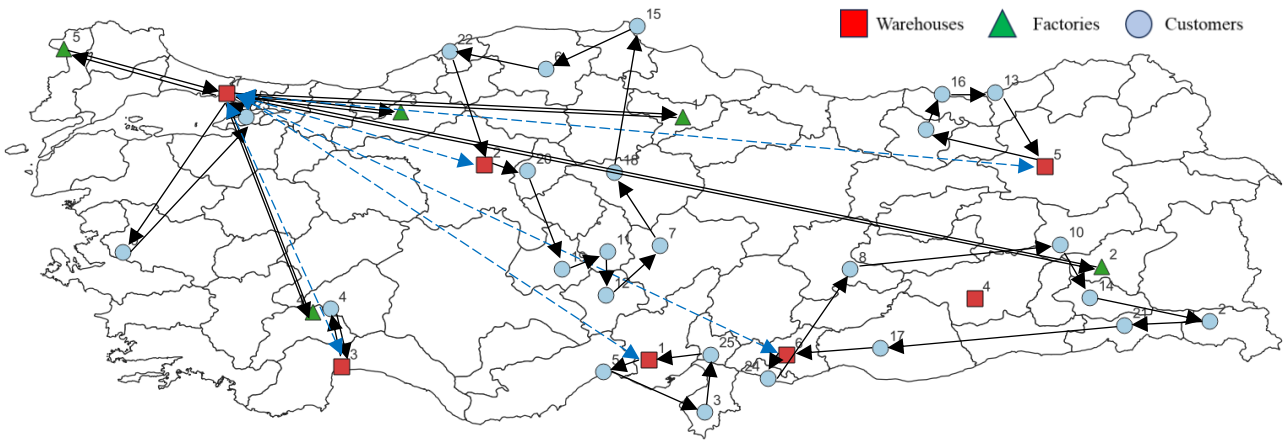

(b) Instance 7-5-25-L

Figure 6. Vehicle routing plots for the best solutions

Table 4. Unit variable costs for the candidate warehouses of both example instances

| Warehouse 1 | Warehouse 2 | Warehouse 3 | Warehouse 4 | Warehouse 5 | Warehouse 6 | Warehouse 7 |
| --- | --- | --- | --- | --- | --- | --- |
| 0.218 | 0.234 | 0.289 | 0.282 | 0.287 | 0.292 | 0.183 |

On the other hand, Figure 7 (a) and (b) shows the values of $\%S_{WI}$ under different warehouse capacities with different warehouse candidate numbers. Figure 7 reveals that for same-sized instances, $\%S_{WI}$ decreases as the warehouse capacities become larger. This conclusion is valid for all tested instances, which suggests that our proposed WSPSDP model performs better when the warehouse capacities are under 'Compact' conditions.

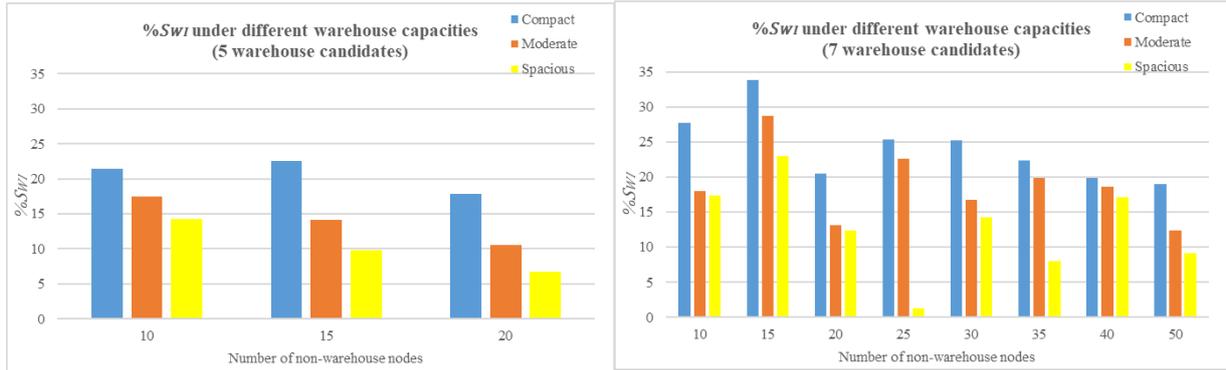

Figure 7. $\%S_{WI}$ values under different warehouse capacities

(a) Instances with 5 warehouse candidates (b) Instances with 7 warehouse candidates

## 6. CONCLUSIONS

In this paper, we proposed a new network for the WSPS that incorporated many new features, such as inter-warehouse transportation and multi-allocation scheme for the collection and delivery processes between warehouses and non-warehouse nodes. Then, we introduced the WSPSDP as a special case of classical MAHLRP and formulated it as a mixed-integer programming problem. An ALNDS heuristic was used to solve the WSPSDP, and the model and algorithm were verified on different test instances.

Based on the results of the computational experiments, we found that for all tested instances, our WSPSDP model outperforms the existing WSPS-WI and WSPS-SA models, proving the proposed model's efficiency and usefulness. In addition, CPLEX is only feasible to solve instances with up to 20 non-warehouse nodes. However, the ALNDS heuristic can find optimal or better solutions for the WSPSDP model within shorter CPU times than CPLEX. Moreover, the results of $\%R_{SD}$ reveal that for medium-and-large-sized instances, the ALNDS heuristic is stable for providing high-quality solutions. In addition, through investigating the effects of warehouse capacities on the solutions, we conclude that our proposed WSPSDP model performs better when the warehouse capacities are under 'Compact' conditions. Furthermore, we found that when the warehouse capacities are sufficiently large, it is possible to use the warehouse with the lowest variable cost for collection processes from factories and the closest warehouses for delivery processes to customers. Hence, the proposed WSPSDP model has the potential to be applied to promote the usage of existing cheap idle warehouses. Future research directions are to consider time aspects, such as time windows, or to expand the current model to incorporate $CO_2$ emissions, thus formulating a bi-objective approach.